\documentclass[a4paper,12pt]{article}
\usepackage{amsmath}
\usepackage{amsthm}
\usepackage{amssymb}
\usepackage{amscd}
\usepackage[dvips]{graphicx}
\newtheorem{definition}{Definition}[section]
\newtheorem{theorem}[definition]{Theorem}
\newtheorem{proposition}[definition]{Proposition}
\newtheorem{lemma}[definition]{Lemma}

\newtheorem{example}[definition]{Example}

\newcommand{\s}{^{\sharp}}
\newcommand{\n}{^{\natural}}

\newcommand{\infgal}{\mathrm{Inf}\text{-}\mathrm{gal}\, }

\newcommand{\com}{\mathbb {C}}

\newcommand{\N}{\mathbb{N}}
\newcommand{\G}{\mathbb{G}}

\newcommand{\Hom}{\mbox{ $\mathrm{Hom}$ }}
\newcommand{\Id}{\mathrm{Id}}

\newcommand{\Q}{\mathbb{Q}}
\newcommand{\spec}{\mathrm{Spec \, }}

\newcommand{\K}{\mathcal{K}}
\newcommand{\eL}{\mathcal{L}}

\newcommand{\gH}{\mathfrak{H}}
\newcommand{\mapue}[1]{%
     \smash{\mathop{%
      \hbox to 1cm{\rightarrowfill}}\limits^{#1}}}

\newcommand{\NCA}{\text{$(NCAlg/L^{\natural})$}}
\begin{document}
\title{Can we quantize Galois theory ?}
\author{Katsunori Saito and 
Hiroshi Umemura \\  
Graduate School of Mathematics \\
Nagoya University\\ \ 
\\
Email\quad {\small 
m07026e@math.nagoya-u.ac.jp and 
 umemura@math.nagoya-u.ac.jp} }
\date{}
\maketitle
\begin{abstract} 
Heiderlich \cite{hei10} tempts us to explore a quantum Galois theory, a Galois theory 
in which quantum group plays the role of algebraic group. 
We give an example of $q$ skew $\sigma$-iterative differential
field  extension whose Galois 
group is a non-commutative and non-co-comutative
Hopf algebra.
\end{abstract} 
\section{Principle of Galois theory}
Heiderich succeeded in combining our general differential Galois theory and Hopf Galois theory for linear differential or difference equations, and by this he suggests a way that would lead us to quantum Galois theory. We give a simple yet non-tivial example in which quantization of Galois theory really happens. This note is a simplified version of \cite{saiume12}. So we refer the readers who are willing to understand details, to the article \cite{saiume12}. 

Theory in science consists of observations and their description. 
Once the method of observation is fixed, it has limits or there are unobservables by the methods. 
The unobservable are hidden in the ambiguity of the method. 
Finer is the method, fewer are the unobservables. 
We might conclude that the finer is the better. 
Since any method has its ambiguity, this does not seem correct. 

According to his testament written the last night before the duel, Galois says that his theory is a theory of ambiguity. 
Maybe he wanted
to assert the following. 
\begin{theorem}[Galois]
When we fix a method of observation and observe an object, then the set of unobservables has a structure. 
We encounter there very often a principal homogeneous space. 
\end{theorem}
\begin{example}\normalfont
We consider a field extension $\Q(i)/\Q$, where $i$ is a complex number satisfying $i^{2} + 1 = 0$. 
Our method of observation is through polynomials $p(x) \in \Q[x]$ with coefficients in $\Q$. 

Let us set $I = - i$ so that $I^{2} + 1 = 0$. Since we have a $\Q$-isomorphism $\varphi \colon \Q[i] \overset{\sim}{\rightarrow}\Q[I] $ of fields such that $\varphi(i) = I$, for every polynomial $p(x) \in \Q[x]$, we have 
\[
\varphi(p(i)) = p(I), 
\]
we can not distinguish $i$ and $I$ by this method. 
In other words, $i$ and $-i = I$ are ambiguous with respect to
this observation method. 
\end{example}
\section{differential Galois theory}
Algebraic differential equation or a differential ring describes a dynamical system over an algebraic variety. 
\begin{example}\normalfont\label{1.21c}
Let $y$ be a variable over $\com$. 
We consider a differential ring $(\com[y],\, \partial)$, where $\partial \colon \com[y] \rightarrow \com[y]$ is a $\com$-derivation such that $\partial(y) = y$ so that $\partial= y d/dy$. 
Formal integration of the vector field $\partial= y d/dy$ gives 
\[
y \mapsto y + \partial (y) X + \frac{1}{2!} \partial^{2}(y) X^{2} + \cdots = \exp(\partial X)y = y \exp X. 
\]
In other words the differential algebra $(\com[y],\, \partial)$ describes the dynamical system
\begin{equation}\label{2.1}
y \mapsto y \exp X 
\end{equation}
on the affine line $\spec \com [y]$, $X$ being the time. 
We observe the dynamical system by algebraic differential equation on the affine line $\spec \com[y]$. Namely our method of observation is algebraic differential equations on $\spec \com[y]$ and the object to be observed is the dynamical system \eqref{2.1}. Let us set 
\begin{equation}\label{2.2}
Y = y \exp X
\end{equation}
for simplicity. $Y$ satisfies the differential equation
\begin{equation}\label{2.3}
y \frac{\partial Y}{\partial y} = Y
\end{equation}
The solution to differential equation to \eqref{2.3} is 
\begin{equation}\label{2.4}
Y = cy,
\end{equation}
$c \in \com$ being a constant. 

So dynamical system \eqref{2.1} observed through algebraic differential equations is nothing but \eqref{2.4};
\begin{equation}\label{2.5}
y \mapsto cy
\end{equation}
If we look at two solutions of the form of \eqref{2.4}
\[
y \mapsto c_{1}y, \quad y\mapsto c_{2}y
\]
The composite is $y \mapsto c_{1}c_{2}y$. 
So Galois group of $(\com[y],\, \partial)/\com$ is the multiplicative group $\mathbb{G}_{m\, \com}$. 
\end{example}
Our general Galois theory \cite{ume96.1}, \cite{ume96.2} depends on this simple principle. 
For a ordinary differential field extension $(L,\, \partial)/(K,\, \partial)$, we attach a kind of group. 
The observation method is the algebraic differential equation on the algebraic variety $V$ over $K$ such that $K(V) = L$, i.e. $V$ is a model of the field extension $L/K$. 
The object to be observed is the dynamical system on $V$ that the differential field extension $L/K$ determines. 
As the set of ambiguity of this observation, we single out the Galois group of the differential field extension $L/K$.
So that the set of unobservables is a principal homogeneous space of the Galois group. 
\section{$q$-analogue}
Since the 19th century, we looked for $q$-analogue of special function. 
\begin{example}\normalfont
We know the hypergeometric series
\begin{equation}\label{3.1}
F(a,\, b,\, c \, ;\, z) = 1 + \frac{a \cdot b}{c \cdot 1}z + \frac{a (a +1)b (b +1)}{c (c +1)\cdot 1 \cdot 2}z^{2} + \cdots
\end{equation}
that is a power series in $z$ parameterized by $a,\, b,\, c \in \com$. We assume $c \neq 0,\, -1,\, -2, \cdots$. 
The hypergeometric series satisfies the hypergeometric differential equation, which is a linear differential equation of the second order over $\com(z)$; 
\begin{equation}\label{3.2}
z(1-z) \frac{d^{2} y}{dz^{2}} - \{( a + b + 1 )z -\gamma \}\frac{dy}{dz} -a b y = 0. 
\end{equation}
To define a $q$-analogue of the hypergeometric series \eqref{3.1}, let $q \neq 0,\, \pm 1$ be a complex number. Heine introduced a series
\[
\varphi(\alpha,\, \beta,\, \gamma\, ; \, z) = 1 + \frac{[\alpha]^{\ast}_{q}[\beta]^{\ast}_{q}}{[\gamma]^{\ast}_{q}[1]^{\ast}_q}z + \frac{[\alpha]^{\ast}_{q}[\alpha + 1]^{\ast}_{q}[\beta]^{\ast}_{q}[\beta + 1]^{\ast}_{q}}{[\gamma]^{\ast}_{q}[\gamma + 1]^{\ast}_{q}[1]^{\ast}_{q}[2]^{\ast}_{q}}z^2 + \cdots , 
\]
where $[\varepsilon]^{\ast}_{q} = (q^{\varepsilon}- q^{-\varepsilon})/(q - q^{-1})$ for a complex number $\varepsilon$. 
We introduce a shift operator 
$Tf(x) = f(qx)$ and a $q$-difference operator $[\theta]^{\ast}_{q} = (T- T^{-1})/(q - q^{-1})$. 
We also introduce $[\theta + \alpha ]^{\ast}_{q}=(q^{\alpha}T- q^{-\alpha}T^{-1})/(q - q^{-1})$ for a complex number $\alpha$. 
The Heine series satisfies difference equation
\begin{equation}\label{3.3}
(x^{-1}[\theta]^{\ast}_{q}[\theta + \gamma -1]^{\ast}_{q} -[\theta + \alpha]^{\ast}_{q}[\theta+\beta]^{\ast}_{q})\varphi = 0. 
\end{equation}
If we notice 
\[
\lim_{q \rightarrow 1} [\varepsilon]^{\ast}_{q} = \varepsilon
\]
When $q \rightarrow 1$, the Heine series
\[
\varphi(\alpha,\, \beta,\, \gamma,\ ; \, z) \rightarrow F(\alpha,\, \beta,\, \gamma,\ ; \, z)
\]
and the difference equation \eqref{3.3} reduces to hypergeometric equation \eqref{3.2}. In this sense the Heine series is a q-analogue of the hypergeometric series. 
According to Galois theory of linear differential equations, for general complex numbers $\alpha, \, \beta,\,  \gamma$, the Galois group of the hypergeometric equation \eqref{3.1} over $\com(z)$ is $\mathrm{GL}_{2}(\com)$ \cite{sinetal03}. 
We have an analogous theory for linear difference equations that tells us for general $\alpha,\, \beta,\, \gamma$, Galois group of the linear difference equation \eqref{3.1} is $\mathrm{GL}_{2}(\com)$. \cite{and01}
\end{example}
The fact that the Heine series is a $q$-analogue of the hypergeometric series means that the Heine series is a deformation of the hypergeometric series parameterized by $q$. It is quite natural to expect that in this procedure Galois group $\mathrm{GL}_{2}(\com)$ would be deformed along $q$ so that we would get a quantization of the algebraic group $\mathrm{GL}_{2}(\com)$. 
Heiderich \cite{hei10} combined our idea with the idea of Takeuchi \cite{amaetal09} and proposed a general Hopf Galois theory 
giving us a hope to discover a Galois theory in which Galois group is a quantum group. 

We show in \cite{ume11} that so far we consider linear functional equations, however non-commutative and non-cocommutative Hopf algebra may be, Galois group is an algebraic group. 

We show by an example of a $q$-skew $\sigma$-differential field extension, of which Galois group is a non-commutative and non-cocommutative Hopf algebra. 
\section{General differential Galois theory}\label{2.2a}

We assume in this section all algebras we consider are $\com$-algebras so that in particular we work in characteristic $0$. 
A differential algebra $(R,\, \delta)$ consists of a ring $R$ and a derivation $\delta \colon R \rightarrow R$. 
So $\delta(a + b) = \delta(a) + \delta(b)$ and $\delta(ab) = \delta(a)b + a\delta(b) $ for any $a,\, b \in R$. 
We often denote the differential algebra $(R,\, \delta)$ simply by $R$, when there is no danger of confusion of the choice of the operator $\delta$. We also have to consider the abstract ring $R$, which we denote by $R\n$. Hence the most strict notation is $R = (R\n,\, \delta)$. 

For a differential algebra $(R,\, \delta)$, we have the universal Taylor morphism
\[
\iota \colon (R,\, \delta) \rightarrow \left( R\s[[X]],\, \frac{d}{dX}\right)
\]
sending an element $y\in R$ to its formal Taylor expansion 
\[
\iota(y) = \sum_{n=0}^{\infty} \frac{1}{n!} \delta^{n}(y)X^{n}, 
\]
where $X$ is a variable over the ring $R\n$. 
The universal Taylor morphism is a differential algebra morphism. 
The universal Taylor morphism defines a map 
\[
\spec R\n[[X]] \rightarrow \spec R,
\]
in other words an $R\n[[X]]$-valued point of $\spec R\n$ tangent to the vector field $\delta$ on $\spec R$ whose initial condition at $X = 0$ is the identity map $\Id_{R} \colon \spec R \rightarrow \spec R$. 
So the universal Taylor morphism describes the formal integration of the vector field $\delta$ on $\spec R$ with the initial condition $\Id \colon \spec R\n \rightarrow \spec R\n$ which is an $R\n$-valued point of $\spec R\n$. See the Example below. 
\begin{example}\normalfont \label{3.4}
Let $y$ be a variable over $\com$ so that $\com[y]$ is the polynomial ring of one variable. 
We consider a differential ring $(\com[y],\, \delta)$ with $\delta = y d/dy$. 
The universal Taylor morphism 
\[
\iota \colon \com[y] \rightarrow \com[y][[X]]
\]
maps the variable $y$ to
\[
\iota(y) = \sum_{n=0}^{\infty}\frac{1}{n!} \delta^{n}(y) X^{n} = y\exp X. 
\]
So $\tilde{Y}(X):= \iota(y) = y \exp X$ satisfies
\begin{equation}\label{3.5}
\frac{d}{dX} Y = Y
\end{equation}
Hence $\tilde{Y}(X)$ is the integration of \eqref{3.5} taking the initial condition $\tilde{Y}(0)=y$. 
\end{example}
We attach to a differential field extension $(L,\, \delta)/(K,\, \delta)$ such that the abstract field $L\n$ is finitely generated over the abstract field $K\n$ Galois group $\infgal(L/K)$. 
Our theory is relative so that the base field $(K,\, \delta)$ is arbitrary but to understand this note the reader may assume $(K,\, \delta) = (\com,\, \delta)$ with $\delta = 0$. 
We have the universal Taylor morphism 
\[
\iota \colon L \rightarrow L\n[[X]]
\]
that is the formal integration of the vector field $\delta$ on the model $V_{K}$ of the field extension $L\n/K\n$. 
To be more precise, we take a model $V_{K}$ of the field extension $L\n/K\n$ so that $V_{K}$ is a $K\n$-scheme of finite type such that $K^{\natural}(V_{K}) \simeq L\n$. 
Then $\delta \colon L \rightarrow L$ defines a rational vector field on the model $V_{K}$ and the universal Taylor morphism defines the formal integration of the vector field on $V_{K}$. 

{\textit We observe the dynamical system defined by the Taylor morphism on $V_{K}$ through algebraic differential equation on the algebraic variety $V_{K}$}. 

To this end, we take a mutually commutative basis $\{D_{1},\, D_{2},\, \cdots ,\, D_{d}\}$ of derivations of the $L\n$-vector space $\mathrm{Der}\,(L\n/K\n)$ of $L\n$ over $K\n$ so that
\begin{enumerate}
\renewcommand{\labelenumi}{(\arabic{enumi})}
\item $d = tr.d. [L\n : K\n]$, 
\item $[D_{i},\, D_{j}] := D_{i}D_{j} - D_{i}D_{j} = 0$ for $1 \leq i,\, j \leq d$
\end{enumerate}
Let us introduce the partial differential field
\[
L\s := (L\n,\, \{D_{1},\, D_{2},\, \cdots ,\, D_{d}\}). 
\]
We add the set of derivations $\{D_{1},\, D_{2},\, \cdots ,\, D_{d}\}$ to the differential ring $(L\n[[X]],\, d/dx)$ such that the $D_{i}$'s operate on the power series ring $L\n[[X]]$ through coefficients. So we get a partial differential ring
\[
L\s [[X]]:= (L\n[[X]],\, \{d/dX,\, D_{1},\, D_{2},\, \cdots ,\, D_{d}\})
\]
\begin{definition}\label{1.21b}
The Galois hull $\eL/\K$ of the given differential field extension is the partial differential algebra extension $\K \subset \eL$ in the partial differential ring $L\s[[X]]$ such that the partial differential algebra $\eL$ is the partial differential algebra of $L\s[[X]]$ generated by $\iota(L)$ and $L\s$ with derivations $d/dX,\, D_{1},\, D_{2},\, \cdots , \, D_{d}$ and such that the partial differential subalgebra $\K$ is differential generated by $i(K)$ and $L\s$ in $L\s[[X]]$. 
\end{definition}
{\bf Continuation of Example \ref{1.21c}. } Let us see what happens in the particular example. 
We take $L = (\com(y),\, yd/dy),\, K = (\com,\, yd/dy) = (\com,\, 0)$ so that $L/K$ is a differential field extension. 
We take as a basis of the $\com(y)$-vector space $\mathrm{Der}\,(\com(y)/\com)$ the derivation $D_{1} =d/dy$ as a basis of one dimensional $\com(y)$-vector space $\mathrm{Der}\,(\com(y)/\com)$. So $yd/dy$ operates on the power series ring $\com(y)[[X]]$ through coefficients so that we may write the generator $D_{1} = \partial/\partial y$. We work in the partial differential ring 
\[
\left( L\s (y) [[X]], \left\{ \frac{\partial}{\partial y},\, \frac{d}{dX} \right\}\right). 
\]
It follows from Definition \ref{1.21b}
that the partial differential ring $\eL$ is generates by $\iota(L) = \com(y\exp X)$ and $\com(y)$ in 
\[
\left(\com(y)[[X]],\, \left\{ \frac{\partial}{\partial y},\, \frac{d}{dX} \right\}\right). 
\] 
Therefore $\eL = \com(y).\com(y\exp X)$ and $\K = \com(y)$. 
Now we are ready to measure the {\it ambiguity} of the inclusion map $\iota \colon \eL \rightarrow L\s[[X]]$ caused of our observation method. 
In the partial case, we similarly define the universal Taylor morphism. 
\begin{equation}\label{1.17a} 
\begin{array}{c}
L\s = (L,\,\{ D_{1},\, D_{2},\, \cdots ,\, D_{d}\})\hfill\\ 
\qquad \rightarrow \left( L\n[[W_{1},\,W_{2},\, \cdots,W_{d}]] \left\{ \dfrac{\partial}{\partial W_{1}} ,\, \dfrac{\partial}{\partial W_{2}}, \cdots, \dfrac{\partial}{\partial W_{d}} \right\} \right)
\end{array}
\end{equation}
where $W_{1},\, W_{2},\cdots,\, W_{d}$ are variables. 
Thanks to the morphism \eqref{1.17a}, we get a morphism
\begin{equation}\label{1.17b}
L\s[[X]] \rightarrow L\n[[W_{1},\, W_{2},\, \cdots,\, W_{d}]][[X]], 
\end{equation}
which is compatible with 
\[
\left \{\frac{d}{dX},\, D_{1},\,D_{2},\, \cdots,\, D_{d}\right\} \, \text{ and }\, \left\{\frac{\partial}{\partial X},\, \frac{\partial}{\partial W_{1}},\, \frac{\partial}{\partial W_{2}},\, \cdots,\, \frac{\partial}{ \partial W_{d}}\right\}. 
\]

Restricting the partial differential morphism \eqref{1.17b} to the partial differential subalgebra $(\eL,\, \{d/dX,\, D_{1},\, D_{2},\, \cdots,\, D_{d} \})$, we get the canonical partial differential morphism
\[ 
\iota \colon \eL \rightarrow L\n[[W_{1},\, W_{2},\, \cdots,\, W_{d}]][[X]] = L\n[[W,\, X]]. 
\]
For an $L\n$-algebra $A$ so that we have a structural morphism
\begin{equation}\label{1.17bc}
L\n \rightarrow A, 
\end{equation}
we get the canonical morphism
\[
\iota \colon \eL \rightarrow L\n[[W,\, X]] \rightarrow A[[W,\, X]]
\]
which we denote also $\iota$. 
\begin{definition}\label{1.17u}
We denote the category of commutative $L\n$ algebras by 
\[
(CAlg/L\n)\]
and the category of sets by $(Set)$. We define a functor
\[
\mathcal{F}_{L/K} \colon (CAlg/L^{\natural}) \rightarrow (Sets)
\]
by setting 
\begin{align*}
&\mathcal{F}_{L/K}(A)=\{ \varphi \colon \eL \rightarrow A[[W,\,X]] \\
&\qquad \qquad \left| \begin{array}{l}
\text{$\varphi$ is a partial differential algebra morphism such that }\\
\text{(1) $\varphi \equiv \iota$ modulo nilpotent elements of $A[[W, \,X ]]$, }\\
\text{(2) $\varphi|_{\K} = \iota|_{\K}$ }
\end{array} \right\}
\end{align*}
for every commutative $L\n$-algebra $A$. 
\end{definition}
The condition (1) in Definition \ref{1.17u} means that for every element $y \in \eL$, the coefficients of the power series $\varphi(y) - \iota(y) \in A[[W,\, X]]$
in the $W_{i}$'s and $X$ are nilpotent element of the ring $A$. 

We define the Galois group $\infgal(L/K)$ as follows. 
\begin{definition}\label{1.21.5}
We define a group functor 
\[
\infgal(L/K) \colon (CAlg/L^{\natural}) \rightarrow (Grp)
\]
by setting
\begin{align*}
&\infgal(L/K)(A) = \{f \colon \eL \hat{\otimes}_{L^{\sharp}} A[[W]] \rightarrow \eL \hat{\otimes}_{L^{\sharp}} A[[W]] \\
& \qquad \qquad \left| \begin{array}{l}
\text{$f$ is an automorphism of $\eL \hat{\otimes}_{L^{\sharp}} A[[W]]$ continuous with} \\
\text{respect to the $W$-adic topology, $f \equiv \Id$ modulo nilpotent }\\
\text{elements and $f$ fixes every element of $\K \hat{\otimes}_{L^{\sharp}} A[[W]]$. }
\end{array}\right\}
\end{align*}
where $\eL \hat{\otimes}_{L^{\sharp}} A[[W]] = \eL \hat{\otimes}_{L^{\sharp}} A[[W_{1},\, W_{2},\, \cdots,\, W_{d}]]$ is the completion of 
\[
\eL \otimes_{L^{\sharp}} A[[W]] = \eL \otimes_{L^{\sharp}} A[[W_{1},\, W_{2},\, \cdots,\, W_{d}]]
\]
with respect to $W = (W_{1},\, W_{2},\, \cdots,\, W_{d})$-adic topology. We use the similar notation for $\K \hat{\otimes}_{L^{\sharp}} A[[W]]$. 
\end{definition}
\begin{theorem}[Umemura\cite{ume96.2}]
The infinitesimal Galois group $\infgal(L/K)$ operates on the functor $\mathcal{F}_{L/K}$ in such a way that $(\infgal(L/K),\, \mathcal{F}_{L/K})$ is a principal homogeneous space. 
\end{theorem}
\section{General difference Galois theory}\label{2.2b}
Let $(R,\, \sigma)$ be a difference ring so that $\sigma \colon R \rightarrow R$ is an endomorphism of commutative $\Q$-algebra $R$. 
A difference ring arises mostly in the following setting. 
Let $S$ be a commutative $\Q$-algebra and $\N = \{0,\, 1,\, 2, \, \cdots\}$ be the set of non-negative integers. 
We denote by
\[
F(\N,\, S):= \{ f \colon \N \rightarrow S\}
\]
the set of functions on $\N$ taking values in the ring $S$. So $F(\N,\, S)$ is a ring. 
It is convenient to express a function in $F(\N,\, S)$ by a matrix
\[
\begin{bmatrix}
0 &1 & 2 & \cdots\\
f(0) & f(1) & f(2) & \cdots
\end{bmatrix}. 
\]
We have the shift operator 
\[
\Sigma \colon F(\N,\, S) \rightarrow F(\N,\, S)
\]
by setting for $f \in F(\N,\, S)$
\[
(\Sigma f)(n) = f(n + 1) \quad \text{ for } n \in \N. 
\]
Therefore $\Sigma$ is an algebra morphism and $(F(\N,\, S),\, \Sigma)$ is a difference ring. 
The difference analogue of the universal Taylor morphism is the universal Euler morphism. 
For a difference ring $(R,\, \sigma)$, the universal Euler morphism is defined by 
\begin{align*}
\iota \colon (R,\, \sigma) &\rightarrow (F(\N,\, R), \, \Sigma)\\
a &\mapsto \begin{bmatrix}
0 & 1 & 2 & \cdots\\
a & \sigma(a) & \sigma^{2}(a) & \cdots
\end{bmatrix}
\end{align*}
that is a difference algebra morphism so that compatible with difference operation $\sigma$ and $\Sigma$. 
General difference Galois theory is established by Shuji 
Morikawa\cite{mori09}. 
The theory depends on the same idea as the general differential theory. Given a difference field extension $(L,\,\sigma)/(K,\, \sigma)$, we simply replace the universal Taylor morphism by universal Euler morphism and apply the idea of \S \ref{2.2a} to construct the Galois hull $\eL/\K$, the infinitesimal deformation functor $\mathcal{F}_{L/K}$ and the Galois group $\infgal(L/K)$. 
\section{Hopf Galois Theory}\label{2.2c}
The universal Taylor morphism and the universal Euler morphism play the prominent role in Galois theory. What are these morphisms? Where do they come from? 

The specialists in Hopf algebra tried to generalize the Picard-Vessiot theory (\cite{swe69}, \cite{amaetal09}). One of their main goals seems to be a uniform understanding of the differential Picard-Vessiot theory and the difference Picard-Vessiot theory. 

The coordinete ring $\com[\G_{a}] = \com[t]$ of the additive algebraic group $\G_{a}$ is a Hopf algebra. 
First of all, it is a polynomial algebra $\com[t]$ of one variable. It has also a co-algebra structure
\[
\Delta \colon \com[\G_{a}] = \com[t] \rightarrow \com[t]\otimes_{\com}\com[t], \; t \mapsto t\otimes 1 + 1 \otimes t 
\]
arising from the group low 
\[
\G_{a} \times \G_{a} \rightarrow \G_{a}. 
\]
We have an algebra morphism $u\colon \com \rightarrow \com[\G_{a}] = \com[t]$ of inclusion and the $\com$-algebra morphism $\varepsilon \colon \com[t]=\com[\G_{a}]\rightarrow \com$ sending $t$ to $0$. 
Moreover we have a $\com$-algebra isomorphism $S \colon \com[\G_{a}]=\com[t]\rightarrow \com[\G_{a}]=\com[t]$ sending $t$ to $-t$. 
So $\com[\G_{a}]$ is a Hopf algebra. 

Let us recall very briefly bi-algebra and Hopf algebra. 
For details, see \cite{swe69}. 
A $\com$-bialgebra $(H,\, m,\, u,\, \Delta,\, \varepsilon)$ has both algebra structure $m \colon H \otimes_{\com} H \rightarrow H,\, a\otimes b \mapsto ab$ and co-algebra structure $\Delta \colon H \rightarrow H\otimes_{\com} H$ so that the maps $m$ and $\Delta$ are $\com$-linear. 
The $\com$-linear map $u \colon \com \rightarrow H$ is the unit morphism and the $\com$-linear map $\varepsilon \colon H \rightarrow \com$ is the co-unit morphism. 
They satisfy a several commutative diagrams. 
A $\com$-Hopf algebra $(H,\, m,\, u,\, \Delta,\, \varepsilon,\, S)$ is a $\com$-bialgebra $(H,\, m,\, u,\, \Delta,\, \varepsilon)$ equipped with a $\com$-linear isomorphism
\[
S : H \rightarrow H
\]
called the antipode satisfying some conditions \cite{swe69}. 
\begin{example}\normalfont
The coordinate ring $\com[G]$ of an affine algebraic group over $\com$ is a commutative Hopf algebra. The universal enveloping algebra $U(\mathfrak{g})$ of a $\com$-Lie algebra $\mathfrak{g}$ is a Hopf algebra that is co-commutative. 
\end{example}
For $\com$-vector space $V,\, W$, We denote the set $\Hom (V,\, W)$ of $\com$-linear maps by ${}_{\com}M(V,\,W)$. It is well-known in linear algebra that we have
\begin{equation}\label{1.25a}
{}_{\com}M(U\otimes V,\, W)\simeq {}_{\com}M(U,\, {}_{\com}M(V,\, W))
\end{equation}
for vector spaces $U,\, V,\, W$. 
We apply \eqref{1.25a} to the particular case where $U=W$ is a $\com$-algebra $R$ and $V$ is a $\com$-Hopf algebra;
\begin{equation}\label{1.25b}
{}_{\com}M(R\otimes_{\com} H,\, R)\simeq {}_{\com}M(R,\, {}_{\com}M(H,\, R))
\end{equation}
We notice here that the dual ${}_{\com}M(H,\, R)$ of the Hopf algebra $H$ is an $R$-algebra. 
\begin{proposition}\label{1.26a}
For a algebra $R$ and the Hopf algebra $H = \com[\G_{a}] = \com[t]$ and a linear map $\Psi \in {}_{\com}M(R\otimes_{\com} H,\, R)$ the following conditions (1),\, (2),\, (3) are equivalent. 
\begin{enumerate}
\renewcommand{\labelenumi}{(\arabic{enumi})}
\item The linear map $\Psi \colon R \otimes_{\com} H \rightarrow R$ is an operator of algebra $H = \com[\G_{a}]$ satisfying
\[
\Psi((ab)h)=\sum_{\alpha \in I} \Psi(a\otimes h_{\alpha})\Psi(b\otimes k_{\alpha}), 
\]
where 
\[
\Delta(h) = \sigma_{\alpha \in I}h_{\alpha}\otimes  k_{\alpha}
\]
with $h_{\alpha},\, k_{\alpha} \in H$, $I$ being a finite index set. 
\item $R$ is a differential algebra with
\[
\delta(a) = \Psi(a\otimes t)
\]
fore every $a \in R$. 
\item The linear map $\Psi' \in {}_{\com}M(R,\, {}_{\com}M(H,\,R))$ corresponding to $\Psi \in {}_{\com}M(R\otimes H,,\, R)$ defines a $\com$-algebra morphism
\[
\Psi' \colon R \rightarrow {}_{\com}M(H,\, R) = R[[X]]. 
\]
\end{enumerate}
\end{proposition}
For the Hopf algebra $\com[\G_{a}]$, the dual algebra ${}_{\com}M(H,\, R)$ is isomorphic to the ring $R[[X]]$ of power series with coefficients in $R$ where the condition (2) is satisfied, the morphism $\Psi' \colon R \rightarrow R[[X]]$ is the universal Taylor morphism. So Proposition \ref{1.26a} characterized the universal Taylor morphism. 
We can generalize Proposition \ref{1.26a} for a bialgebra $H$. 
\begin{proposition}\label{1.26b}
Let $H$ be a $\com$-algebra. Except for this, we keep the notation in Proposition \ref{1.26a}. 
For a $\com$-algebra $R$ and a linear map $\Psi \in {}_{\com}M(R\otimes_{\com}H,\, R)$, the following conditions are equivalent. 
\begin{enumerate}
\renewcommand{\labelenumi}{(\arabic{enumi})}
\item The assertion (1) in Proposition \ref{1.26a} holds. 
\item $\Psi' \colon R \rightarrow {}_{\com}M(H,\, R)$ is a $\com$-algebra morphism. 
\end{enumerate}
\end{proposition}
If we take the $\com$-Hopf algebra $H = \com[\G_{a}]$, then the condition (1) in Proposition \ref{1.26b} says $(R,\, \delta)$ is a difference algebra and the condition (2) gives the universal Euler morphism (\cite{hei10})
\section{$q$-skew iterative $\sigma$-differential algebra}
We study in section \ref{2.2a} differential theory and in section \ref{2.2b} difference theory. 
We explaind the idea of unifying both cases for linear equations in section \ref{2.2c}. Heiderich \cite{hei10} pointed out that we can apply the idea of Hopf algebraists in Picard-Vessiot theory to our
Galois theory for non-linear equations. 
A new landscape is open when we deal with neither commutative nor co-commutative Hopf algebras. As an example of such Hopf algebra, we consider $q$-skew iterative $\sigma$-algebras. For the hopf algebra in question, see \cite{hei10}. 
For a positive integer $n$ and an element $q\neq 0,\, 1$ of a commutative ring $C$, we set
\[
[n]_{q} := \sum_{i=0}^{n-1}q^{i} \in C
\]
that is equal to
\[
\frac{q^{n} - 1}{q-1}
\]
if $C$ is a field. It is convenient to set
\[
[0]_{q}=1
\]
For a non-negative integer, we introduce the $q$-factorial
\[
[n]_{q}! = [n]_{q}[n-1]_{q} \cdots [1]_{q} 
\]
and 
\[
[0]_{q}!=1. 
\]
The $q$-binomial coefficient is defined by 
\[
\binom{i + j}{i}_{q} = \binom{i + j}{j}_{q} = \frac{[i+j]_{q}!}{[i]_{q}![j]_{q}!}
\]
for non-negative integers $i,\,j$. We can show that
\[\binom{i + j}{i}_{q} = \binom{i + j}{j}_{q} \]
is a well-determines element in the ring $C$. See Kac and Cheung \cite{kacetal02}. 
\begin{definition}\label{1.26d} Let $q$ be a non-zero complex number. 
A $q$-skew iterative $\sigma$-differential algebra ($qsi$ algebra for short) consists of an associative $\com$-algebra $R$ that is not necessarily commutative, a $\com$-algebra endomorphism $\sigma \colon R \rightarrow R$ and a set $\{ \theta^{(i)}\}_{i\in \N}$ of $\com$-linear maps
\[
\theta^{(i)} \colon R \rightarrow R\quad \text{ for } i = 0,\, 1,\, 2,\, \cdots
\]
satisfying the following conditions. 
\begin{enumerate}
\renewcommand{\labelenumi}{(\arabic{enumi})}
\item $ \theta ^{(0)} = \Id_{R}$. 
\item $  q^i \sigma \theta ^{(i)} = \theta ^{(i)}\sigma \quad \text{ for every } i = 0,\, 1,\, 2, \cdots$. 
\item $ \displaystyle{\theta ^{(l)}(ab) =\sum_{l = m+n , \, m, n \in \N} \sigma^{m}(\theta^{(l)}(a) ) \theta^{(m)}(b)}$,
\item  $\theta^{(i)}\circ \theta^{(j)} =\binom{i+j}{i}_{q} \theta^{(i+j)} \quad \text{ for } m,\, n \in \N$. 
\end{enumerate}
\end{definition}
\begin{example}\normalfont \label{1.30a}
A difference ring $(R,\,\sigma$ has a trivial $qsi$ algebra structure. In fact, it is sufficient to define $\theta^{(0)} = \Id_{R}$ and $\theta^{(i)} = 0$ for $i \geq 1$. Then $R,\, \sigma, \, \theta^{\ast}$ is a $qsi$ algebra for every $0 \neq q \in \com$. 
\end{example}
\begin{example}\normalfont
A differential algebra $(R,\, \delta)$ in characteristic $0$ also has a trivial $qsi$ algebra structure for $q = 1$. We set $\theta^{0} = \Id_{R}$ and $\theta^{(i)} = (1/i!)\delta^{i}$ for $i \geq 1$. 
Then$(R,\, \Id,\, \theta^{\ast})$ is a $qsi$ algebra for $q = 1$. 
\end{example}
These example shows that $qsi$ algebra generalizes both difference and differential algebra.
\begin{example}\normalfont \label{2.6a}
Non-canonical $qsi$ algebra structures on a difference ring $(R,\, \sigma)$. We choose an element $\lambda \in R$. We set
\begin{align*}
\theta^{(0)} = \Id_{R}, \quad \theta^{(1)} = \lambda (\sigma - \Id_{R}), \\
\theta^{(i)} = \frac{1}{[i]_{q}!}(\theta^{(1)})^{i} \quad \text{ for } i=1,\,2,\, \cdots,
\end{align*}
so that the $\theta^{(i)}$'s are linear maps $R\rightarrow R$. 
If the element $\lambda \in R$ satisfies $q\sigma(\lambda) = \lambda$, then $(R,\, \sigma,\, \theta^{\ast})$ is a $qsi$ algebra. 

In particular if the difference field $(\com(t),\, \sigma)$ is a difference subalgebra of $(R,\, \sigma)$ such that $\sigma(t) = qt$ with $q \neq 1,\, 0$, then we can take
\[
\lambda = \frac{1}{(q-1)t}
\]
so that
\[
\theta^{(1)} = \frac{1}{(q-1)t}(\sigma - \Id_{R}). 
\]
\end{example}
\begin{example}\normalfont
[Non-commtative $qsi$ algebra] Let $(R,\, \sigma)$ be a difference algebr that may be non-commutative. Let $q$ be a non-zero complex number. 
$(R,\, \sigma)[[X]]$ is a twisted power series ring so that it is a power series ring $R[[X]]$ as an additive group with the following commutation relations;
\[
aX = X\sigma(a) \quad \text{ for every } a \in R
\]
so that
\[
a X^{n} = x^{n}\sigma^{n}(a) \quad \text{ for } n=1,\, 2,\,3 \cdots. 
\]
Hence the ring $(R,\, \sigma)[[X]]$ is non-commutative even if $R$ is a commutative ring. 
We notice that every element of $(R,\, \sigma)[[X]]$ is uniquely written in the form of power series $ \sum_{i=0}^{\infty} X^{i}a_{i}$ with $a_{i} \in R$ for $i \in \N$. 
We show that $(R,\, \sigma)[[X]]$ has a structure of $qsi$ algebra. 
The ring endomorphism $\Sigma \colon (R,\, \sigma)[[X]] \rightarrow (R,\, \sigma)[[X]]$ is given by
\[
\Sigma \left( \sum_{i=0}^{\infty} X^{i}a_{i} \right) = \sum_{i=0}^{\infty} X^{i}q^{i}\sigma(a_{i})
\]
for an element $\sum_{i=0}^{\infty}X^{i}a_{i} \in R[[X]]$. 
The definition of the linear map 
\[
\Theta^{(i)}\colon (R,\, \sigma)[[X]]  \rightarrow (R,\, \sigma)[[X]]\]
is given by
\[
\Theta^{(1)}\left( \sum_{i=0}^{\infty} X^{i}a_{i} \right) = \sum_{i=0}^{\infty}[i]_{q}X^{i-1}a_{i}
\]
for $\sum_{i=0}^{\infty} X^{i}a_{i} \in (R\, \sigma)[[X]]$. 
\[
\Theta^{(0)} = \Id_{(R,\, \sigma)[[X]]} \text{ and } \Theta^{(i)} = \frac{1}{[i]_{q}!} \Theta^{i} \text{ for } i = 2,\, 3,\, \cdots. \]
Then $((R,\, \sigma)[[X]], \, \sigma,\, \Theta^{\ast})$ is a $asi$ algebra that is non-commutative in general. 
\end{example}
A $qsi$ algebra is characterized by a Hopf algebra as a differential ring is characterized by the Hopf algebra $\com[\G_{a}]$. See Heiderich \cite{hei10}, Proposition 2.3.15. So we can replace the universal Taylor morphism by the universal Hopf morphism for a $qsi$ algebra. 
\begin{definition}[Heiderich\cite{hei10}, Proposition 2.3.17]
Let $(R,\, \sigma,\, \Theta^{\ast})$ is a $qsi$ algebra. The universal Hopf morphism
\[
\iota \colon (R,\, \sigma,\, \theta^{\ast}) \rightarrow (F(\N,\, R)[[X]],\, \Sigma,\, \Theta^{\ast})
\]
is a morphism of $qsi$ algebras given by
\begin{equation}\label{2.6b}
\iota(a) = \sum_{i=0}^{\infty} X^{i}u[\theta^{(i)}(a)], 
\end{equation}
for every element $a \in R$, 
where for an element $f \in R$, we denote the function
\[
\begin{bmatrix}
0 & 1 & 2 & \cdots\\
f & \sigma(f) & \sigma^{2}(f) & \cdots
\end{bmatrix} \in F(\N,/, R)
\]
by $u[f]$. 
\end{definition}
The universal Hopf morphism generalizes both the universal Taylor morphism and Euler morphism. 
\section{Example}\label{sec8}
\begin{example}\label{8.1}\normalfont
For a complex number $q \neq 0$, we consider a $\com$-difference field $(\com(t),\, \sigma)$, where $\sigma \colon \com(t) \rightarrow \com(t)$ is a $\com$-automorphism of the rational function field $\com(t)$ of one variable such that 
\begin{equation} \label{2.4a}
\sigma(t) = qt.
\end{equation}
So $(\com(t),\, \sigma)/(\com,\, \Id_{\com})$ is a difference field extension. 
\eqref{2.4a} shows that the extension $(\com(t),\, \sigma)/(\com,\, \Id_{\com})$ is a Picard-Vessiot extension and the Galois group $\mathrm{Gal}((\com(t),\, \sigma)/(\com,\, \Id_{\com}))$ that is the group of $\com$-difference automorphism of $(\com(t),\, \sigma)$, is isomorphic to $\G_{m,\,\com}$. 
Let us introduce trivial derivations $\delta^{\ast}$ on $(\com(t),\, \sigma)$ as in Example \ref{1.30a}. 
So 
\[
\delta^{(0)} = \Id_{\com(t)},\quad \delta^{(i)} = 0 \quad \text{ for } i=1,\, 2,\, \cdots. 
\]
We get $(\com(t),\, \sigma,\, \delta^{\ast})$ is a $qsi$ field. 
The $qsi$ field extension
\[
(\com(t),\, \sigma,\, \delta^{\ast})/(\com,\, \Id_{\com},\, \delta^{\ast})
\]
is a Picard-Vessiot extension and its Galois group is isomorphic to $\G_{m,\, \com}$. See Heiderich \cite{hei10},\, Umemura \cite{ume11}. 

We now introduce $\theta^{\ast}$ as in Example \ref{2.6a}. 
\[
\theta^{0} = \Id_{\com(t)},\quad \theta^{(1)} (a) = \frac{\sigma(a) -a}{(q-1)t} \quad \text{ for } a \in \com(t). 
\]
In general, we set
\[
\theta^{(i)} = \frac{1}{[i]_{q}!}(\theta^{(1)})^{i} \quad \text{ for } i = 1,\, 2,\, \cdots
\]
Therefore $(\com(t),\, \sigma,\, \theta^{\ast})$ is a $qsi$ field and the subfield $\com$ of $\com(t)$ consisting of constant polynomials is invariant under $\sigma$ and $\theta^{\ast}$. 
So $\com$ is a $qsi$ subfield of $(\com(t),\, \sigma,\, \theta^{\ast})$. Let us study the $qsi$ field extension $L/K:= (\com(t),\, \sigma,\, \theta^{\ast})/\com$. We will see later the extension $L/K$ is not a Picard-Vessiot extension. 
Intuitively, the generator $t$ of $\com(t)$ over $\com$ satisfies the following system of difference-differential equations
\[
\sigma(t) = qt,\quad \theta^{(0)}(t) = t,\quad \theta^{1}(t) = 1,\quad \theta^(i)(t) = 0\quad\text{ for } i \geq 2. 
\]
that is not linear homogeneous. 
We show that the Galois group $\mathrm{Gal}(L/K)$ is a Hopf algebra that is neither commutative nor
co-commutative. According to \cite{hei10}, we know that if $L/K$ were a Picard-Vessiot extension, the Galois group $\mathrm{Gal}(L/K)$ would be a commutative Hopf algebra. 

Let us first write the universal Hopf morphism
\[
\iota \colon L \rightarrow F(\N,\, L\n)[[X]]
\]
so that 
\[
\iota(a) = \sum_{n=0}^{\infty}x^{n}u[\theta^{n}(a)] \quad \text{ for every } a \in L
\]
as we learned in \eqref{2.6b}. 
In particular, we have
\begin{align}
\iota(t) &= u[t] + X[u(1)] \notag \\
&=
\begin{bmatrix}
0 & 1 & 2 & \cdots\\
t & qt & q^{2}t & \cdots
\end{bmatrix}
+ X\begin{bmatrix}
0 & 1 & 2 & \cdots\\
1 & 1 & 1 & \cdots
\end{bmatrix}\notag \\
&= t
\begin{bmatrix}
0 & 1 & 2 & \cdots\\
1 & q & q^{2} & \cdots
\end{bmatrix}
+ X \notag \\
&= tQ + X, \label{2.6c}
\end{align}
if we introduce the function
\[
Q := \begin{bmatrix}
0 & 1 & 2 & \cdots \\
1 & q & q^{2} & \cdots
\end{bmatrix}
\in F(\N,\, L\n) \subset F(\N,\, L\n)[[X]]
\]
We identify subrings through inclusions
\[
L\n \subset F(\N,\, L\n) \subset F(\N,\, L\n)[[X]]. 
\]
We take the derivation $d/dt \in \mathrm{Der}\,(L\n/K\n)$ as a basis of $1$-dimensional $L\n$-vector space $\mathrm{Der}\,(L\n/K\n)$. 
So $L\s = (L\n,\, d/dt)$. 
Therefore $\eL$ is generated by $\iota(L)$ and $L\s$ invariant under the operators $\Sigma,\, \Theta^{\ast}, d/dt$ of the ring $F(\N,\, L\s)[[X]]$. 
It follows from \eqref{2.6c}
\[
\frac{d}{dt}\iota(t) = \frac{d}{dt}(tQ + X) = Q. 
\]
So $Q\in \eL$ and $\iota(t) -tQ = X \in \eL$. 
Therefore $\eL \supset L\s \langle Q,\, X \rangle_{g}$, the subalgebra generated by $Q$ and $X$ over $L\s$ in $F(\N,\, L\s)[[X]]$. 
To be more precise 
\[
\eL = L\s ( tQ + X )\langle X,\, Q \rangle_{g}. 
\]
We consider the universal Taylor morphism
\[
L\s = (L\n,\, d/dt) \rightarrow (L\s[[W]],\, d/dW), 
\]
which gives us
\[
\eL \rightarrow F(\N,\, L\s)[[X]] \rightarrow F\left(\N,\, \left(L\n[[W]],\, \frac{d}{dW} \right) \right)[[X]], 
\]
which we denote by $\iota$. 
Since $\eL$ is non-commutative, it is more natural to extend functors $\mathcal{F}_{L/K}$ and $\infgal(L/K)$ on the category $(CAlg/L\n)$ of commutative $L\n$-algebras to the category of non-commutatibe algebrs. 
To this end, we introduce the category \NCA \, of not necessarily commutative $L\n$-algebra such that $L\n$ is contained in the center of the algebra. 
For an algebra $A \in ob(NCAlg/L\n)$, we consider an infinitesimal deformation
\[
h \colon \eL \rightarrow F(\N,\, L\n[[W]])[[X]] \rightarrow F(\N,\, A[[W]])[[X]]
\]
of the canonical morphism $\iota$. The deformation is determined by the images $h(Q)$ and $h(X)$. 
Since $\Sigma Q = q Q,\; dQ/dt = 0$ and $\Theta^{(1)}Q=0 $ we have
\[
h(Q) = eQ, 
\]
where $e \in A$ such that $e-1 \in A$ nilpotent. 
Similarly as $dX/dt=0,\; \Sigma X = qX,\; \Theta^{(1)} X = 1$, 
\[
h(X) = fQ + X,
\]
where $f \in A$ is nilpotent. 
Since
\[
XQq = X\Sigma (Q) = QX, 
\]
we must have
\[
h(XQq) = h(QX)
\]
so that we have
\[
h(X)h(Q)q = h(Q)h(X),
\]
namely
\[
(fQ+X)eQq = eQ(tQ+X). 
\]
We must have
\[
feq = eq
\]
conversely we can take any $e,\, f \in A$ such that $e-1$ and $f$ are nilpotent element of $A$ and such that$ef=feq$. 
\end{example}
\section{Quantum group $\mathfrak{h}_{q}$}
We encountered a quantum group in section \ref{sec8}. 
\begin{definition}
We work in the category $(NCAlg/\mathbb{C})$. Let $A$ be a not necessarily commutative $\mathbb{C}$-algebra. 
We say that  two subsets $S,\,  T$ of $A$ are mutually commutative if for every $s\in T,\, 
 t \in T$, we have $[s,t]=st-ts = 0$.
\end{definition}
We can prove the following two Lemmas by direct calculations. 
\begin{lemma} For two matrices
\[
Z\sb{1} = \begin{bmatrix}
e\sb{1} &f\sb{1} \\
 0 & 1 
\end{bmatrix},\quad  
Z\sb{2} = \begin{bmatrix}
e\sb{2} & f\sb{2}\\
0       & 1 
\end{bmatrix} \in H\sb{q}(A), 
\]
if $\{e\sb{1},\, f\sb{1}\}$ and $\{e\sb{2},\, 
f\sb{2}\}$ are mutually commutative, then the product matrix 
\[
Z\sb{1} Z\sb{2} \in H\sb{q}(A).
\]
\end{lemma}
\begin{lemma}
For a matrix
\[
Z = \begin{bmatrix}
e & f\\
0 & 1
\end{bmatrix} \in H\sb{q}(A),
\]
if we set
\[
\tilde{Z} = \begin{bmatrix}
e^{-1} & -e^{-1}f \\
 0  & 1
\end{bmatrix} \in M\sb{2},
\]
then 
\[
\tilde{Z} \in H\sb{q^{-1}}(A) \text{ and } 
\tilde{Z}Z = Z\tilde{Z} = I\sb{2}. 
 \]
\end{lemma}
According to a general principle of Manin \cite{man88} section 1, we can find a $\com$-Hopf algebra $\mathfrak{h}_{q}$ as follows. 
\begin{lemma}
Let $u$ and $v$ be symbols over $\com$.
We have shown that we find a $\mathbb{C}$-Hopf algebra 
\begin{equation} \label{hopfa1}
\mathfrak{h}\sb q = \mathbb{C}\langle u,\, u^{-1}, \, v\rangle \sb{alg}/(uv -q^{-1}vu)
\end{equation}
as an algebra so that
$$
uu^{-1} = u^{-1} u = 1 , \qquad u^{-1} v = qvu^{-1}.
$$
Definition of the algebra $\gH\sb q$ determines the multiplication
\[
m \colon \gH \sb q\otimes\sb{\mathbb{C}}\gH\sb q \to \gH\sb q, 
\]
the unit 
\[
\eta \colon \mathbb{C} \to \gH\sb q, 
\]
that is the composition of natural
 morphisms $$\com \rightarrow \com\langle u, \, u^{-1}, \, v\rangle\sb{alg}$$ and $$\com\langle u, \, u^{-1}, \, v \rangle\sb{alg} \rightarrow \com\langle u, \, u^{-1}, \, v\rangle\sb{alg}/(uv-q^{-1}vu)=\gH\sb q.$$ 
 The product of matrices gives 
the co-multiplication
\[
\Delta  \colon \gH\sb q \to \gH\sb q\otimes\sb{\mathbb{C}}\gH\sb q, 
\] 
that is a $\com$-algebra morphism 
defined by
\[
\Delta (u)=u\otimes u, \qquad \Delta (u^{-1})=
u^{-1}\otimes u^{-1}, 
\qquad \Delta (v)=u\otimes v + v\otimes 1,
\]
for the generators $u,\, u^{-1},  \, v$ of the algebra ${\gH\sb q}$, 
the co-unit is a $\com$-algebra morphism 
\[
\epsilon \colon \gH\sb q \to \mathbb{C}, \qquad \epsilon(u)=\epsilon (u^{-1})=1, \, \epsilon(v)=0 
 \]
 for the generators $u, \, u^{-1}, \, v$ of the algebra 
${\gH\sb q}$.
The antipode 
$$
i\colon {\gH\sb q} \to {\gH\sb q}
$$
is a $\com$-anti-algebra morphism given by 
$$
i(u)=u^{-1}, \qquad i(u^{-1})= u,  \qquad i(v) = -u^{-1}v. 
$$
\end{lemma}
\begin{proposition}\label{9.1a}
So we have for $A \in ob(NCAlg/\com)$, the set $\gH_{q}(A)$ of $A$-valued points of the quantum group $\gH_{q}$ that is identified with
\[
H_{q}(A) = \left\{ \left.
\begin{bmatrix}
e & f \\
0 & 1
\end{bmatrix} \right| e, f\in A,\, qef=fe\,\text{ and } e \, \text{ is invertible }
\right\}. 
\]
\end{proposition}
\section{Conclusion}
We conclude from section \ref{sec8}, Example \ref{8.1} and Proposition \ref{9.1a}, that the  formal completion $\gH_{q} \hat{\otimes}_{\com} L\n$ of the quantum group at the unit element, which is considered to be a functor $(NCAlg/L\n) \rightarrow (Set)$ operates on the functor $\mathcal{F}_{L/K} \colon (NCAlg/L\n) \rightarrow (Set)$ in such a way that we can regard 
\[
(\gH\hat{\otimes}_{\com}L\n,\, \mathcal{F}_{L/K})
\]
is a principal homogeneous space. 
The Galois group of the $qsi$ field extension $\com(t)/\com$ is the quantum group $\gH\hat{\otimes}_{\com}L\n$. 
\bibliographystyle{plain.bst}
\bibliography{umemura2b}
\end{document}